\documentclass[a4paper,12pt]{scrartcl} 
\title{Covariate balance under Bayesian decision theory}

\author{
    André A. F. Fumis \hspace{1mm} \texttt{fumis@ime.usp.br}, \\
    Victor Fossaluza \hspace{1mm} \texttt{victorf@ime.usp.br}, \\
    Rafael B. Stern \hspace{1mm} \texttt{rbstern@gmail.com} \\
    Institute of Mathematics and Statistics, University of São Paulo
}
\date{September 1, 2025}

\usepackage[T1]{fontenc}
\usepackage{nameref}
\usepackage{algorithm, algorithmic}
\usepackage[colorlinks=true, urlcolor=blue, citecolor=blue, linkcolor=blue]{hyperref}

\usepackage{amsfonts, amsmath, amssymb, amsthm, thmtools, mathtools} 
\usepackage{mathrsfs}
\usepackage{cleveref}
\usepackage{accents}
\usepackage{nicematrix}

\theoremstyle{plain}

\newtheorem{theorem}{Theorem}[section]
\newtheorem{lemma}[theorem]{Lemma}
\newtheorem{corollary}[theorem]{Corollary}
\theoremstyle{definition}

\newtheorem{specification}[theorem]{Model Specification}

\crefname{section}{section}{sections}
\Crefname{section}{Section}{Sections}
\crefname{table}{table}{tables}
\Crefname{table}{Table}{Tables}

\usepackage{color, enumitem, float, layout, multirow, setspace, verbatim} 
\setlist[enumerate]{leftmargin=*}
\usepackage[colorinlistoftodos]{todonotes}
\usepackage{listings}
\usepackage{caption}
\usepackage{epigraph}
\usepackage{subfig}
\usepackage{ulem}
\usepackage{makecell}

\usepackage[round]{natbib}

\usepackage{graphicx}
\graphicspath{{figures/}{../figures/}}
\newcommand{\X}{\mathbf{X}}
\newcommand{\Xc}{\mathbf{X_C}}
\newcommand{\Xt}{\mathbf{X_T}}
\newcommand{\U}{\mathbf{U}}
\newcommand{\Uc}{\mathbf{U_C}}
\newcommand{\Ut}{\mathbf{U_T}}
\newcommand{\Z}{\mathbf{Z}}
\newcommand{\J}{\mathbf{J}}
\newcommand{\G}{\mathbf{G}}
\newcommand{\Oo}{\mathbf{O}}
\newcommand{\w}{\mathbf{w}}

\newcommand{\y}{\mathbf{y}}
\newcommand{\R}{\mathbf{R}}
\newcommand{\A}{\mathbf{A}}
\newcommand{\B}{\mathbf{B}}
\newcommand{\bvec}{\mathbf{b}}
\newcommand{\D}{\mathbf{D}}
\newcommand{\Hh}{\mathbf{H}}
\newcommand{\Q}{\mathbf{Q}}
\newcommand{\V}{\mathbf{V}}
\newcommand{\Sca}{\mathbf{S}}
\newcommand{\zerovec}{\mathbf{\undertilde{0}}}
\newcommand{\onevec}{\mathbf{\undertilde{1}}}
\newcommand{\xb}{\mathbf{\bar{x}}}
\newcommand{\xbc}{\mathbf{\bar{x}_C}}
\newcommand{\xbt}{\mathbf{\bar{x}_T}}

\newcommand{\gb}{\mathbf{\bar{g}}}
\newcommand{\gbc}{\mathbf{\bar{g}_C}}
\newcommand{\gbt}{\mathbf{\bar{g}_T}}

\newcommand{\obc}{\mathbf{\bar{o}_C}}
\newcommand{\obt}{\mathbf{\bar{o}_T}}
\newcommand{\uu}{\mathbf{u}}
\newcommand{\K}{\mathbf{K}}
\newcommand{\Ll}{\boldsymbol{\Lambda}}
\newcommand{\indsim}{\stackrel{ind.}{\sim}}

\begin{document}

\maketitle

\begin{abstract}%
We study optimal sample allocation between treatment and control groups under Bayesian linear models. We derive an analytic expression for the Bayes risk, which depends jointly on sample size and covariate mean balance across groups. Under a flat conditional prior, the covariate mean balance term simplifies to the Mahalanobis distance. Our results reveal that the optimal allocation does not always correspond to equal sample sizes, and we provide sufficient conditions under which equal allocation is optimal. Finally, we extend the analysis to sequential settings with groups of patients arriving over time.
\end{abstract}


\section{Introduction}
\label{sec:intro}

Controlled experiments such as clinical trials are essential to modern science as a straightforward method of estimating the causal effect of an intervention. The simplest design consists of splitting the experimental units into two groups, treatment and control, and applying distinct interventions to each group. If the allocation is performed completely at random, the expected value of the difference between groups is the average treatment effect \citep[ATE;][]{rubin_estimating_1974}. It is true even in the presence of important confounding variables, as the random allocation balances them on average. However, any particular allocation may still exhibit imbalances on confounding variables. Considering experiments are often performed only a single time before evaluation, a need arises to minimize all sources of variation, allowing for precise inference on the effect of interest.

A common form of controlling the variation is blocking. While excellent in situations with a large amount of control over the experimental units and only a modest amount of categorical variables to be balanced, the approach falters when dealing with many covariates due to the curse of dimensionality. Also, it forces inappropriate dicotomization for continuous variables. Considering such problems, newer techniques propose diminishing a proper choice of distance or other dissimilarity criterion between the treatment and control groups, where the criterion is based on the variables being balanced.

One situation extensively discussed in the literature is single allocations, when the decision maker has a set amount of experimental units to divide into treatment and control, whose covariate data is all available at the moment of allocation. In that situation, the Mahalanobis distance between covariate means has taken place as the most common dissimilarity criterion under linear model assumptions \citep{cox_use_1957, cox_randomization_1982, morgan_rerandomization_2012, morgan_rerandomization_2015, lauretto_haphazard_2017, kallus_optimal_2018}.

\citet{cox_use_1957, cox_randomization_1982} evaluates the average variance of the analysis of covariance estimator, which relies on a linear adjustment on the covariates, over the randomization distribution of treatment assignments for equal sample sizes. The work exposes the connection between covariance imbalance and the estimator variance. An initial discussion is provided by \citet{cox_randomization_1982} on the validity of inference based on the randomization distribution when the realized design may be atypical, outlining the basis for inference based on rerandomization.

\citet{morgan_rerandomization_2012, morgan_rerandomization_2015} propose the allocation of treatment assignments based on rerandomization of treatment assignments with equally sized groups until the Mahalanobis distance between covariate means is smaller than a threshold, based on the theoretical quantiles of the distance statistic over the randomization distribution. Using the difference-in-means estimator, they show a reduction in variance of the estimator when the model fits a linear structure. Asymptotic properties of allocation based on rerandomization are given by \citet{li_asymptotic_2018} and \citet{schultzberg_asymptotic_2021}.

In both the randomization distribution of \citet{cox_use_1957, cox_randomization_1982} and the rerandomization distribution of \citet{morgan_rerandomization_2012, morgan_rerandomization_2015}, an important aspect for inference is the existence of a non-degenerate distribution of estimates from a fixed sample by leveraging different potential allocations, which makes frequentist confidence intervals and null hypothesis significance tests possible. Another approach, haphazard intentional allocation \citep{lauretto_haphazard_2017}, performs strict optimization of a convex combination of the Mahalanobis distance applied to the data and to random noise, while providing computational techniques for generating sufficiently balanced allocations in viable time. The addition of noise serves the same purpose of (re)randomization instead of strict optimization in previous works, inserting randomness to the allocation to avoid degenerate distributions of estimates.

\citet{kallus_optimal_2018} likewise considers the difference-in-means estimator for balanced group sizes, studying the estimator variance over the randomization distribution with fixed covariates and potential outcomes. The approach taken is that of minimax decision theory, considering adversarial covariate and outcome values. Under a linear structure for the conditional outcome mean, \citet{kallus_optimal_2018} shows the minimax decision is to fully minimize the Mahalanobis distance. Other allocation methods are shown to be minimax optimal under different model structure assumptions, including complete randomization when no assumptions are made.

Our work, like others, studies sample allocation under linear model assumptions. We take the approach of Bayesian decision theory, obtaining an analytic formula for the decision risk of an allocation, given by expected parameter variance. The formula is valid for an interpretable conjugate class of prior distributions of model hyperparameters, allowing the use of domain knowledge when allocating the sample. In the flat conditional prior limit, our formula recovers the Mahalanobis distance between covariate means as a key component of the decision risk. Furthermore, the formula makes explicit the components of decision risk attributable to covariate imbalance, through the Mahalanobis distance, and to treatment group size differences.

Bayesian decision theory for sample allocation has seen comparatively little discussion next to frequentist approaches. \citet{kasy_why_2016} studies optimal allocation under square loss for a wide class of models, showing the computational tractability of loss evaluation for a proposed allocation and the independence of potential outcomes and treatment assignment conditional to covariates used for the allocation, an important result for identification of causal effects. Our work follows the same framework, but focusing on the interpretation of the prior, the attribution of treatment group size differences and covariate imbalance to the risk and the characterization of the optimal allocation, while making parametric assumptions.

The optimization and rerandomization approaches have the limitation of requiring all experimental units to be available for allocation at once, which is unrealistic for many experiments, such as clinical trials that depend on patients developing the condition to be treated. Existing sequential approaches to covariate balancing employ stochastic algorithms for greedy reduction of the dissimilarity criterion between groups \citep{pocock_sequential_1975, stern_combining_2015, qin_adaptive_2022}. We likewise extend our sample allocation results to a greedy sequential approach, showing the formula for decision risk obtained for single allocations can be adapted to context of the allocation of a subsample after the rest of the sample has already been allocated.

The use of Bayesian decision theory leads to non-randomized, systematic designs, a topic which has been very discussed. As mentioned before, randomization plays an important role in frequentist inference. In pure Bayesian theory, however, the sampling design, including any randomization involved, plays no part in the analysis after the sample has been collected. The sampling design is relevant for minimizing the expected loss function of the optimal estimator chosen for the analysis, but any two estimates will be identical if the same sample is collected from different sampling designs. \citet{solomon_optimal_1970} give a review of a related problem to sample allocation, optimal sample selection, from the perspectives of both minimax and Bayesian decision theory. \citet{kasy_why_2016} advocates for deterministic designs in sample allocation even when considering minimax risk, under the assumption model parameters may be adversarial to assignment mechanism but not to the specific allocation chosen. \citet{banerjee_decision_2017} unifies minimax and Bayesian decision theory by considering optimal design of experiments relative to a combination of the expected loss against the designer's prior and the largest expected loss against a set of different, possibly adversarial, priors, reflecting the subjective opinions of an audience of fellow stakeholders. Notably, the choice of experiment, ranging from the designer's preferred one to a fully randomized one, depends on the available sample size, as external validity becomes more important relative to the gain in precision from optimal allocation relative to full randomization. We follow \citet{kasy_why_2016} in considering the optimal, deterministic design, unlike some designs mentioned \citep{morgan_rerandomization_2012, lauretto_haphazard_2017, cox_use_1957, cox_randomization_1982} that avoid strict optimization or intentionally introduce stochasticity in order to perform valid frequentist inference. Our proposed Bayesian model leverages the independence of parameters and assignment mechanism of treatments, conditional on observed data, to enable inference regardless of experiment design.

The paper is organized as follows: In Section \ref{sec:homo-model-statement} we state a Bayesian model for experiment outcomes with respect to observed covariates, and study optimal decision making for sample allocation under the model. We present two main theorems, at the ends of Subsections \ref{subsec:general-prior} and \ref{subsec:sequential}, containing analytic formulas for the decision risk of an allocation of all experimental units at once and of sequential allocation of remaining experimental units after some have already been allocated. Subsection~\ref{subsec:flat-prior} relates the main theorem of Subsection~\ref{subsec:general-prior} to the existing literature by examining the limiting flat conditional prior case. Subsection~\ref{subsec:group-size} examines the validity of the common assumption of equal group sizes in covariate balance. It presents a sufficient condition under which the Bayes-optimal decision lies within this constraint, and also provides a counterexample demonstrating that the optimal decision may contain unbalanced group sizes.
\section{Covariate balance under the linear homoscedastic model}
\label{sec:homo-model-statement}

In this paper we present the relationship between the decision risk of an allocation of experimental units between treatment and control groups and two characteristics of the allocation, the treatment arm sample sizes and the covariate balance. We begin by defining some notation for the paper, and then proceed to declare a statistical model and give a general expression for the risk of an allocation under the model. In following subsections we will develop this expression for the risk to show the role of the treatment arm sample sizes and covariate balance.

As shown by \citet{kallus_optimal_2018}, model structure plays a crucial part in decision-theoretic optimal sample allocation. We consider the Normal, homoscedastic, linear model structure under the potential outcomes causal framework given by
\begin{equation}
    y_i(j) \indsim N(\gamma_j + \mathbf{x}_i\beta, \sigma^2),\;j=0,1,
\end{equation}

where $y_i(0)$ represents the potential outcome for experimental unit $i$ under the control intervention, $y_i(1)$ represents the potential outcome under the treatment intervention and $y_i = (1-w_i) \cdot y_i(0) + w_i \cdot y_i(1)$, where $w_i \in \{0,1\}$ is the treatment indicator. This model admits a matrix representation
\begin{equation}
\label{likelihood}
    \y \sim N(\Z \zeta, \sigma^2 \mathbf{I}_n),
\end{equation}
where $\Z = \begin{pmatrix}
    \mathbf{1}-\w & \w & \X
\end{pmatrix}$ and $\zeta = \begin{pmatrix}
    \gamma_0 \\
    \gamma_1 \\
    \beta
\end{pmatrix}$.

The quantity of interest of a controlled experiment is often the average treatment effect $\tau = E[Y(1) - Y(0)] = \gamma_1 - \gamma_0 = \A \zeta$, where $\A = \begin{pmatrix} -1 & 1 & 0 & \dots & 0 \end{pmatrix}$. We wish to find the optimal allocation of experimental units $\w$ between the control and treatment arms. The allocation of units to each experimental arm happens after observation of covariates $x_1, \dots, x_n$, but before observation of outcomes $y_1, \dots, y_n$, as those depend on the choice of intervention.

The allocation $\w$ can be analyzed under Bayesian decision theory. Consider a loss function given by the squared error of the optimal estimate under allocation $\w$. Then the risk associated with allocation $\w$ after observation of covariate matrix $\X$ is
\begin{gather}
    \label{risk_definition}
    r_{\X}(\w) = E[Var(\tau \mid \y, \X, \w, \sigma^2) \mid \X, \w] \; \citep[p.~228]{degroot_optimal_2004}.
\end{gather}

Equation (\ref{risk_definition}) presents the risk for the Normal linear model, but treats expectations implicitly with respect to the prior. We restrict our attention to the conjugate prior distribution for regression parameters $\zeta$ and $\sigma^2$ given by $\zeta, \sigma^2 \sim NIG(\zeta_0, \V_0, a_0, b_0)$ such that $E[\sigma^2]$ is finite (satisfied by $a_0>1, b_0>0$). One particular case of interest is the flat conditional prior limit, discussed in Subsection~\ref{subsec:flat-prior}. A precise statement of the full statistical model is
\begin{specification} 
    \label{spec:bayesian_model_spec}
    \;
    \begin{itemize}
        \item $\mathbf{z}_i = (1-w_i, w_i, x_i)$ deterministic function of $\mathbf{x}_i$, $w_i$
        \item $y_i \mid \mathbf{x}_i, w_i, \zeta, \sigma^2 \sim c.i.i.d\; N(\mathbf{z}_i \zeta, \sigma^2)$
        \item $\mathbf{x}_i \mid \psi \sim c.i.i.d\; F_{\psi}$; $F_{\psi}$ arbitrary
        \item Parameters $\zeta, \sigma^2, \psi \sim f(\zeta, \sigma^2, \psi) = f(\zeta, \sigma^2)f(\psi)$, \textit{i.e.}, $(\zeta, \sigma^2)$ and $\psi$ are \textit{a priori} independent
        \item $\zeta, \sigma^2 \sim NIG(\zeta_0, \V_0, a_0, b_0)$ with $a_0>1, b_0>0$.
    \end{itemize}
\end{specification}
Under such conditions, then
$\zeta, \sigma^2 \mid \w, \X, \y \sim NIG(\zeta_1, \V_1, a_1, b_1)$ \cite[p.~252]{degroot_optimal_2004} where:
\begin{itemize} 
	\item $\zeta_1 = \V_1(\V_0^{-1}\zeta_0 + \Z'\Z\hat{\zeta})$
	\item $\V_1 = (\V_0^{-1} + \Z'\Z)^{-1}$
	\item $a_1 = a_0 + \frac{n}{2}$
	\item $b_1 = b_0 + \frac{\zeta_0'\V_0^{-1}\zeta_0 + \y'\y - \zeta_1'\V_1^{-1}\zeta_1}{2}$
	\item and $\hat{\zeta} = (\Z'\Z)^{-1}\Z'\y$ is the traditional OLS estimator.
\end{itemize}

This model formulation allows us to write the decision risk in an analytically tractable manner that depends only on one prior hyperparameter, $\V_0$, as well as on the observed samples $\X$:

\begin{lemma}
    \label{lem:bayesian_risk}
    Under Model Specification~\ref{spec:bayesian_model_spec}, $r_{\X}(\w) = \A \V_1 \A' E[\sigma^2]$.
\end{lemma}

Accordingly, the Bayes decision is the allocation which minimizes $\A \V_1 \A'$ and depends on $\X$ and $\V_0$ only. $\V_0$ is the prior covariance matrix for regression coefficients $\zeta$, up to a constant. The other prior hyperparameters, the mean for $\zeta$ given by $\zeta_0$ and the hyperparameters $a_0$ and $b_0$ that specify $\sigma^2$'s marginal distribution, influence posterior inference but not the allocation decision.
\subsection{General prior}
\label{subsec:general-prior}

The matrix $\V_0$ is confined to the restrictions of a covariance matrix, which are symmetry and positive-semidefiniteness. In particular, it has the following block matrix representation:
\begin{equation} 
    \V_0 =
    \begin{pmatrix}
        \nu & \varrho \\
        \varrho' & \gamma
    \end{pmatrix}
\end{equation}
where $\nu$ and $\gamma$ are symmetric and positive-semidefinite. The $\nu$ block stands for the covariance matrix of the intercepts for the control and group arms, the $\gamma$ block for the covariance matrix of the coefficients associated with covariates in the model and the $\varrho$ block for the covariances between intercepts and covariate associated coefficients.

Suppose $\V_0$ is positive-definite and, equivalently, invertible. Then $\V_0^{-1}$, known as the precision matrix, is also positive-definite and admits a Cholesky decomposition $\V_0^{-1} = \Q'\Q$, where
\begin{equation} 
    \Q =
    \begin{pmatrix}
        \Hh & \B \\
        \mathbf{0} & \D
    \end{pmatrix},
    \label{V0_as_QtQ}
\end{equation}
$\Hh$ is $2 \times 2$, $\B$ is $2 \times p$, $\D$ is $p \times p$ and $\Hh$ and $\D$ are upper-triangular. We make the additional restriction that $\Hh$ must be diagonal, giving it the form $\Hh = \begin{pmatrix}
    h_1 & 0 \\
    0 & h_2
\end{pmatrix}$. Using block matrix inversion \citep[p.~108]{bernstein_matrix_2009}, we obtain
\begin{itemize} 
    \item $\D'\D = \gamma^{-1}$, meaning $\D$ is the Cholesky decomposition of $\gamma^{-1}$;
    \item $\Hh = (\nu - \varrho \gamma^{-1} \varrho')^{-\frac{1}{2}}$, which exists and is diagonal $\iff$ $\nu - \varrho \gamma^{-1} \varrho'$ is diagonal;
    \item $\B = -\Hh \varrho \gamma^{-1}$.
\end{itemize}

Note that $\B$ is a $2 \times p$ matrix, and denote its first row by $\bvec_1$ and the second by $\bvec_2$.

This representation of the prior covariance matrix $\V_0$ allows us to obtain an analytical representation for the risk.

Let $n_C$ be the number of units assigned to the control arm and $n_T$ the number assigned to the treatment arm, $n_C + n_T = n$. Let also $\xbc$ and $\xbt$ be the $1 \times p$ row vectors of the covariate means in the control and treatment arms. Denote the scatter matrix of covariate matrix $\X$ by $\Sca(\X) = \X'\X - n \xb'\xb$, and consider the other following constants:
\begin{itemize} 
    \item $s_C = n_C + h_1^2$;
    \item $s_T = n_T + h_2^2$;
    \item $s=s_C+s_T$;
    \item $\gbc = \frac{h_1 \cdot \bvec_1 + n_C \cdot \xbc}{s_C}$;
    \item $\gbt = \frac{h_2 \cdot \bvec_2 + n_T \cdot \xbt}{s_T}$;
    \item $\gb = \frac{h_1 \cdot \bvec_1 + h_2 \cdot \bvec_2 + n \xb}{s}$.
\end{itemize}

\begin{theorem} 
    \label{thm:risk_theorem_with_prior}
    Under a prior $\zeta, \sigma^2 \sim NIG(\zeta_0, \V_0, a_0, b_0)$ where $\V_0^{-1} = \Q'\Q$ as given above and fixed $n_C$ and $n_T$,
    \begin{gather*}
        r_{\X}(\w) = \frac{(\frac{s}{s_C s_T})^2}{\frac{s}{s_C s_T} - (\gbt - \gbc)(\X'\X + \B'\B - s \gb'\gb + \D'\D)^{-1}(\gbt - \gbc)'} E[\sigma^2].
    \end{gather*}
\end{theorem}

One important aspect of Theorem~\ref{thm:risk_theorem_with_prior} is that the risk is monotonically increasing on a quadratic form whose matrix does not depend on the allocation of units. This means the computational problem of obtaining the optimal allocation is equivalent to that of minimizing the quadratic form, and that evaluation of the quadratic form does not depend on a different matrix inversion per allocation, but on a single inversion per sample.

While the representation seems convoluted, consider the specific case when $\Hh^2$ is composed of integer diagonal values $h_1^2$ and $h_2^2$ and define 
\begin{gather}
    \G = \begin{pmatrix}
        \X \\
        \bvec_1/h_1 \\
        \vdots \\
        \bvec_1/h_1 \\
        \bvec_2/h_2 \\
        \vdots \\
        \bvec_2/h_2
    \end{pmatrix},    
\end{gather}
where row $\bvec_1/h_1$ is repeated $h_1^2$ times and row $\bvec_2/h_2$ is repeated $h_2^2$ times.

\begin{corollary} 
    \label{col:risk_corollary_with_prior}
    Under a prior $\zeta, \sigma^2 \sim NIG(\zeta_0, \V_0, a_0, b_0)$ where $\V_0^{-1} = \Q'\Q$ as given in \eqref{V0_as_QtQ} with $\Hh^2$ composed of integer diagonal values and fixed $n_C$ and $n_T$,
    \begin{gather*}
        r_{\X}(\w) = \frac{(\frac{s}{s_C s_T})^2}{\frac{s}{s_C s_T} - (\gbt - \gbc)(\Sca(\G)+\D'\D)^{-1}(\gbt - \gbc)'} E[\sigma^2].
    \end{gather*}
\end{corollary}

This representation gives a possible interpretation to the submatrices $\Hh$, $\B$ and $\D$ that compose $\Q$, the upper-triangular matrix obtained from the Cholesky decomposition of $\V_0^{-1}$, the prior precison matrix. $\G$, obtained as an extension of sample covariate matrix $\X$ with $h_1^2$ repeated observations of value $\frac{\bvec_1}{h_1}$ and $h_2^2$ repeated observations of value $\frac{\bvec_2}{h_2}$, stands as the \textit{de facto} sample covariate matrix, while $\D'\D$ acts as a regularization term. $\Hh^{-1}\B = \begin{pmatrix}
    \frac{\bvec_1}{h_1} \\
    \frac{\bvec_2}{h_2}
\end{pmatrix}$ thus gives the covariate mean vectors for the pseudo-samples in the control and treatment groups, while $\Hh^2 = \begin{pmatrix}
    h_1^2 & 0 \\
    0 & h_2^2
\end{pmatrix}$ gives the sample sizes.
\subsection{Flat conditional prior}
\label{subsec:flat-prior}

A common choice of prior distribution for applied Bayesian data analysis is of those considered non-informative, under which fall the Jeffreys prior \citep{jeffreys_invariant_1946} and the flat prior. For the linear regression model in Model Specification~\ref{spec:bayesian_model_spec}, the two classes of priors mentioned have flat conditional priors $f(\zeta \mid \sigma^2) \; \alpha \; 1$, but $E[\sigma^2]$ is not finite as required by the assumptions stated in Lemma~\ref{lem:bayesian_risk}. A choice of prior which has finite $E[\sigma^2]$ is the $\zeta, \sigma^2 \sim NIG(\zeta_0, \V_0, a_0, b_0)$ distribution under the flat conditional prior limit $\V_0^{-1} \rightarrow \mathbf{0}_{(p+2) \times (p+2)}$, equivalent to $f(\zeta \mid \sigma^2) \; \alpha \; 1$. Note that, per Lemma~\ref{lem:bayesian_risk}, the allocation decision $\w$ depends on the prior only through $\V_0$.

\begin{corollary} 
    \label{col:risk_flat_prior}
    Under the flat conditional prior limit $\V_0^{-1} \rightarrow \mathbf{0}_{(p+2) \times (p+2)}$ and fixed $n_C$ and $n_T$, 
    \begin{gather*}
        r_{\X}(\w) = \frac{(\frac{n}{n_C n_T})^2}{\frac{n}{n_C n_T} - (\xbt - \xbc)\Sca(\X)^{-1}(\xbt - \xbc)'} E[\sigma^2].
    \end{gather*}
\end{corollary}

As an immediate consequence, we see that the Bayes decision in the flat prior limit agrees with previous methods from the (re)randomization distribution framework. The Mahalanobis distance between covariate means is defined by \citet{morgan_rerandomization_2012} as $M = n \frac{n_C}{n} \frac{n_T}{n} (\xbt - \xbc)Cov(\X)^{-1}(\xbt - \xbc)'$ and used as a scalar measure of multivariate covariate balance for the rerandomization criterion, where $Cov(\X)$ is the sample covariance matrix of $\X$ (with denominator $n-1$).

\begin{corollary} 
    \label{col:risk_mahalanobis}
    Under a flat conditional prior $f(\zeta \mid \sigma^2) \; \alpha \; 1$ and fixed $n_C$ and $n_T$, $r_{\X}(\w)=\frac{n}{n_C n_T} \frac{1}{1-\frac{M}{n-1}} E[\sigma^2]$.
\end{corollary}

\begin{corollary}
    \label{col:lauretto}
    The allocation defined in \citet{lauretto_haphazard_2017} with $\lambda = 0$ is the Bayes decision with respect to the statistical model defined by \eqref{likelihood} and the flat conditional prior for $\zeta$ under the restriction of fixed $n_C$ and $n_T$.
\end{corollary}

\citet{morgan_rerandomization_2012} show conditioning the allocation on having a Mahalanobis distance between covariate means smaller than a threshold provides a variance reduction to the difference-in-means estimator. Corollary~\ref{col:risk_mahalanobis} highlights the relationship between the Mahalanobis distance and the decision risk of an allocation when the estimator of choice is the difference between linear regression coefficients, which is the Bayes estimator under the stated conditions, and provides a finite sample optimality result under the treatment group size restriction. Notably, the same dissimilarity metric arises as the object of limitation or minimization under two different estimators. The equivalence between the Bayes decision in \cite{lauretto_haphazard_2017} and \eqref{risk_definition} allows for the use of the described mixed-integer linear programming techniques to solve the original ATE variance minimization problem.
\subsection{Sequential allocations}
\label{subsec:sequential}

In many experiments, such as clinical trials in health applications, experimental units arrive sequentially and must be allocated to the treatment or control arm at different moments. For practical or ethical reasons, waiting for the full sample may be infeasible. Such a reason may be starting treatment of ill patients as soon as possible for better outcomes, for example. In such situations, practitioners aiming to achieve covariate balance may follow adaptive procedures that take into account the patients already allocated. An instance of such a trial is described in \citet{fossaluza_sequential_2009}, where an algorithm based on Aitchinson's compositional distance was used to allocate patients sequentially in a clinical trial comparing pharmacological and psychotherapeutic interventions in OCD patients. The trial featured the additional difficulty of restrictions on the amount of patients allocated to each treatment arm at each point, due to logistical constraints from the psychotherapeutic treatment cohort starting points.

Consider, after allocation $\w_1$ is performed on $n$ units with covariate matrix $\X$ observing outcomes $\y_1$, a second allocation $\w_2$ is performed on $m$ units with covariate matrix $\U$ observing outcomes $\y_2$. Assume $\U$ is sampled $ciid$ from the same distribution $F_{\psi}$ in Model Specification~\ref{spec:bayesian_model_spec} that $\X$ is sampled from. The goal is to balance the treatment and control groups, under the restriction of the original sample's allocation. As such, the posterior distribution $\zeta, \sigma^2 \mid \w1, \X, \y_1 \sim NIG(\zeta_1, \V_1, a_1, b_1)$ is used as the prior for the second allocation, and thus for the allocation decision $\w_2$. By the conjugacy property of the $NIG$ distribution, $\zeta, \sigma^2 \mid \w_1, \X, \y_1, \w_2, \U, \y_2 \sim NIG(\zeta_2, \V_2, a_2, b_2)$, where $\V_2$ and consequently the decision $\w_2$ depend only on $\V_1$ and $\U$.

Let $\U_C$ be the submatrix of $\U$ composed by the $m_C$ units allocated to control and $\U_T$ the one composed by the $m_T$ allocated to treatment by decision $\w_2$. Define the following quantities from the joint sample composed of $\X$ and $\U$:
\begin{itemize} 
    \item $\Oo = \begin{pNiceArray}{cc|c}
        \Block{2-1}<\Large>{\onevec} & \Block{2-1}<\Large>{\zerovec} & \Xc \\
         & & \Uc \\
         \Block{2-1}<\Large>{\zerovec} & \Block{2-1}<\Large>{\onevec} & \Xt \\
         & & \Ut \\
    \end{pNiceArray}$;
    \item $l_C = n_C + m_C$;
    \item $l_T = n_T + m_T$;
    \item $\obc = \frac{\onevec'\Xc + \onevec'\Uc}{l_C}$;
    \item $\obt = \frac{\onevec'\Xt + \onevec'\Ut}{l_T}$.
\end{itemize}
Let now, analogously to Theorem~\ref{thm:risk_theorem_with_prior},
\begin{itemize} 
    \item $s_C = l_C + h_1^2$;
    \item $s_T = l_T + h_2^2$;
    \item $s=s_C+s_T$;
    \item $\gbc = \frac{h_1 \cdot \bvec_1 + l_C \cdot \obc}{s_C}$;
    \item $\gbt = \frac{h_2 \cdot \bvec_2 + l_T \cdot \obt}{s_T}$;
    \item $\gb = \frac{h_1 \cdot \bvec_1 + h_2 \cdot \bvec_2 + l_C \cdot \obc + l_T \cdot \obt}{s}$.
\end{itemize}

\begin{theorem} 
    \label{thm:sequential_risk}
    Under a prior $\zeta, \sigma^2 \sim NIG(\zeta_0, \V_0, a_0, b_0)$ with $\V_0^{-1} = \Q'\Q$ as given in \eqref{V0_as_QtQ} for the first allocation and under fixed $m_C$ and $m_T$ for the second,
    \begin{gather*}
         r_{\U, \X, \y_1, \w_1}(\w_2) = \frac{(\frac{s}{s_C s_T})^2}{\frac{s}{s_C s_T} - (\gbt - \gbc)(\Oo'\Oo+\B'\B-s\gb'\gb+\D'\D)^{-1}(\gbt - \gbc)'} E[\sigma^2 \mid \X, \y_1, \w_1].
    \end{gather*}
\end{theorem}

Theorem~\ref{thm:sequential_risk} shows that $r_{\U, \X, \y_1, \w_1}(\w_2)$ is minimized by the same partial allocation $\w_2$ which would minimize the risk obtained from a single allocation $\w = \begin{pmatrix}
    \w_1 \\
    \w_2
\end{pmatrix}$ over the joint covariate matrix $\begin{pmatrix}
    \X \\
    \U
\end{pmatrix}$, under a restriction of fixed $\w_1$. Simplifications similar to Corollary~\ref{col:risk_flat_prior} apply. This provides a greedy approach to the covariate balancing problem under sequential allocation of experimental units: once some units have already been allocated and other units join the experiment, the best possible approach is to minimize the conditional risk within the constraints of the already allocated units.
\subsection{Treatment group size}
\label{subsec:group-size}

The optimality of the allocation that minimizes Mahalanobis distance between covariate means stands only with a flat conditional prior and under the fixed treatment group size restriction, when $n_T$ and $n_C=n-n_T$ are fixed. When $n_C$ and $n_T$ are allowed to vary, the decision risk does not have a strictly monotonic relationship with the Mahalanobis distance, as the $\frac{n}{n_C n_T}$ term takes different values. The most intuitive split under homoscedasticity of outcomes is equal division between treatment and control. The intuition is further supported by the fact $\frac{n}{n_C n_T}$ is minimized when $n_C=n_T=\frac{n}{2}$, taking value $\frac{4}{n}$.

A sufficient condition for an even split being optimal has been derived in terms of the projection matrix $\Hh(\X)=\X(\X'\X)^{-1}\X'$, known as the hat matrix in linear regression for projecting outcome vector $\y$ to fitted values $\hat{\y}$, the linear projections of $\y$ in the column space of $\X$.

\begin{theorem}
    \label{thm:equal_split_condition}
    Under a flat conditional prior $f(\zeta \mid \sigma^2) \; \alpha \; 1$ and with even $n$, the optimal allocation will have $n_C=n_T=\frac{n}{2}$ if there exists an allocation vector $\w$ with $n_C=n_T=\frac{n}{2}$ such that $(\w-\frac{1}{2}\onevec)'\Hh(\X)(\w-\frac{1}{2}\onevec) \leq \frac{1}{n}$.
\end{theorem}

Theorem~\ref{thm:equal_split_condition} provides a sufficient condition for the optimal allocation to have an equal amount of units in the treatment and control groups. The condition may not always be reached. In such cases, it is possible that the optimal allocation, as measured by the expected posterior variance for $\tau$, will have a different number of units in the treatment and control groups. Due to the symmetry imposed by the statistical model, notably the homoscedasticity of outcomes, such allocations with the larger amount of units in either the treatment or the control group are equally good. One example of a sample where an equal split is not optimal is Table~\ref{unequal_split_sample}.

\begin{table}[ht]
\centering
\begin{tabular}{rrrr}
  \hline
 $X_1$ & $X_2$ & $X_3$ \\ 
  \hline
0.1 & -0.8 & -1.3 \\ 
0.5 & 2.1 & 1.3 \\ 
0.8 & -0.2 & 0.2 \\ 
-0.3 & 0.3 & 0.6 \\ 
1.1 & -0.8 & 0.0 \\ 
-0.5 & 0.7 & -0.7 \\ 
-0.8 & 1.2 & -0.4 \\ 
-0.7 & 1.0 & 1.4 \\ 
   \hline
\end{tabular}
\caption{A sample where an equal split is not optimal.}
\label{unequal_split_sample}
\end{table}

It was obtained by simulating independent observations from a standard Normal distribution for each of the $8 \cdot 3 = 24$ cells until a sample was found where the optimal allocation did not have an equal split, which took just 21 iterations of a loop in the R software \citep{r_software} under the chosen random seed, and then rounding the numbers to one decimal place for ease of presentation. The optimal allocation splits experimental units into a group with rows 1, 2 and 4 and a group with rows 3, 5, 6, 7 and 8. The smallest amount reached for $(\w-\frac{1}{2}\onevec)'\Hh(\X)(\w-\frac{1}{2}\onevec)$ among equal splits is, rounded to two decimal places, 0.24, while the necessary value to be obtained for the sufficient condition in Theorem~\ref{thm:equal_split_condition} is $\frac{1}{8}=0.125$.
\section{Final remarks}
\label{sec:final}

The methods presented in this paper add to the existing literature \citep{cox_use_1957, cox_randomization_1982, morgan_rerandomization_2012, kallus_optimal_2018} that establishes the Mahalanobis distance as the desired dissimilarity criterion for covariate balance under linearity conditions by proving the variance of the Bayes estimator is a monotonic function of it at the flat prior limit, while generalizing the decision risk for a class of prior distributions. The use of Bayesian inference also has benefits in post-treatment inference, as posterior distributions obtained from the data are independent of the assignment mechanism and may explicitly incorporate domain expertise through the choice of prior distributions. Logically, this comes at the cost of model assumptions, as shown by \citet{kallus_optimal_2018}. Nevertheless, the use of Bayesian inference allows the analyst to choose their own methodology and assumptions irrespective of the designer's choice of sample allocation, while perhaps recognizing the allocation may not have been the optimal choice for the chosen analysis, which likely would've happened as well with a randomized decision.

Strict optimization of the decision risk is a computationally intensive problem for more than a moderate amount of sample units. In practice, practitioners can use techniques such as choosing the best among a number of rerandomizations \citep{morgan_rerandomization_2012} and solving for an approximation of the Mahalanobis distance \citep{lauretto_haphazard_2017}, as a small enough Mahalanobis distance will result in an estimator with near-optimal variance.

Last, we provide a greedy algorithm for allocation of units arriving as sequential groups with the property that, conditionally on past allocations, covariate balance after the current allocation is as best as can be. It provides a guideline to practitioners running experiments where treatment may not be postponed until the full sample is obtained.

A topic for future research is removing the homoscedasticity assumption between the two potential outcomes, allowing for different variances which may be related through a prior distribution, and studying how the treatment group size decision varies with the choice of prior distribution. As shown in Subsection~\ref{subsec:group-size}, in the homoscedastic case there can be particular samples where an equal split is not optimal, but there exists a sufficient condition for the optimality which depends only on the set of allocations with an equal split.

\section*{Acknowledgements and funding}

This study was financed in part by the Coordenação de Aperfeiçoamento de Pessoal de Nível Superior – Brasil (CAPES) – Finance Code 001. The authors are grateful to Marcelo Lauretto for helpful suggestions and code samples.

\bibliography{biblio}

\appendix
\section*{Appendix A - \Cref{sec:homo-model-statement} Proofs}

\begin{proof}[Proof of Lemma~\ref{lem:bayesian_risk}]
    From Model Specification~\ref{spec:bayesian_model_spec}, it follows that \cite[p.~354]{gelman_bayesian_1995}:
    \begin{itemize}
        \item $f(\zeta, \sigma^2, \psi \mid \X, \w, \y) = f(\psi \mid \X)f(\zeta, \sigma^2 \mid \X, \w, \y)$
    	\item $f(\zeta, \sigma^2 \mid \X, \w) = f(\zeta, \sigma^2)$
    	\item $f(\zeta, \sigma^2 \mid \X, \w, \y) \; \alpha \; f(\zeta, \sigma^2)f(\y \mid \X, \w, \zeta, \sigma^2)$.
    \end{itemize}
    
    We obtain
    \begin{align} 
        r_{\X}(\w)
        & = E[Var(\tau \mid \y, \X, \w, \sigma^2) \mid \X, \w]  \nonumber\\
        & = \int_{\sigma^2} \int_{\y} Var(\tau \mid \y, \X, \w, \sigma^2) f(\y \mid \X, \w, \sigma^2) d\y f(\sigma^2 \mid \X, \w) d\sigma^2  \label{risk_1}\\
        & = \int_{\sigma^2} \int_{\y} \A Var(\zeta \mid \y, \X, \w, \sigma^2) \A' \; f(\y \mid \X, \w, \sigma^2) d\y f(\sigma^2) d\sigma^2, \nonumber
    \end{align}
    where a random variable under the integral sign indicates the integration region is the entire support of the random variable.
    
    One property of the $NIG$ distribution is that if $\mathbf{x}, t \sim NIG(\mu, \V, a, b)$, then $\mathbf{x} \mid t \sim N(\mu, \V t)$ \cite[p.~252]{degroot_optimal_2004}. It follows that $\zeta \mid \sigma^2, \w, \X, \y \sim N(\zeta_1, \V_1 \sigma^2)$ and $Var(\zeta \mid \y, \X, \w, \sigma^2) = \V_1 \sigma^2$, where $\V_1$ depends only on the prior, $\X$ and $\w$. Thus, (\ref{risk_1}) simplifies into
    \begin{align} 
        r_{\X}(\w)
        & = \int_{\sigma^2} \int_{\y} \A (\V_1 \sigma^2) \A' \; f(\y \mid \X, \w, \sigma^2) d\y f(\sigma^2) d\sigma^2  \nonumber\\
        & = \A \V_1 \A' \int_{\sigma^2} \int_{\y} f(\y \mid \X, \w, \sigma^2) d\y \; \sigma^2 f(\sigma^2) d\sigma^2 \\
        & = \A \V_1 \A' E[\sigma^2]. \nonumber
    \end{align}
\end{proof}

\begin{proof}[Proof of Theorem~\ref{thm:risk_theorem_with_prior}]
    By Lemma~\ref{lem:bayesian_risk}, our goal is to obtain $\A \V_1 \A'$, where
    $\A = \begin{pmatrix} -1 & 1 & 0 & \dots & 0 \end{pmatrix}$. For that intent, only the upper left $2 \times 2$ block of $\V_1$ is relevant, which can be obtained by block matrix inversion \cite[p.~108]{bernstein_matrix_2009}.

    Without loss of generality, assume $\w = (0, \dots, 0, 1, \dots, 1)'$, such that the first $n_C$ elements of $\w$ are assigned to the control arm and the other $n_T$ are assigned to the treatment arm. We can thus write
    \begin{equation} 
        \Z
        = \begin{pNiceArray}{cc|ccc}
            1 & 0 & \Block{6-3}<\LARGE>{\X} \\
            \vdots & \vdots & & & \\
            1 & 0 & & & \\
            0 & 1 & & & \\
            \vdots & \vdots & & & \\
            0 & 1 & & &
        \end{pNiceArray}
        = \begin{pNiceArray}{cc|c}
            \onevec & \zerovec & \Xc \\
            \zerovec & \onevec & \Xt
        \end{pNiceArray}
    \end{equation}
    where $\zerovec$ and $\onevec$ represent conformable column vectors of 0s and 1s with length implied by context, and $\Xc$ and $\Xt$ are the submatrices defined by the first $n_C$ and last $n_T$ rows of $\X$. Consider the following block matrix representation:
    \begin{equation} 
        \Z'\Z =
        \begin{pNiceArray}{cc}
            \onevec' & \zerovec' \\
            \zerovec' & \onevec' \\
            \hline
            \Xc' & \Xt'
        \end{pNiceArray}
        \begin{pNiceArray}{cc|c}
            \onevec & \zerovec & \Xc \\
            \zerovec & \onevec & \Xt
        \end{pNiceArray}
        = \begin{pNiceArray}{cc|c}
            n_C & 0 & \onevec'\Xc \\
            0 & n_T & \onevec'\Xt \\
            \hline
            \Xc'\onevec & \Xt'\onevec & \Xc'\Xc + \Xt'\Xt
        \end{pNiceArray}.
    \end{equation}

    Together with Expression (\ref{V0_as_QtQ}) of the main paper, we obtain
    \begin{align} 
        \V_1
        & = (\V_0^{-1} + \Z'\Z)^{-1} = (\Q'\Q + \Z'\Z)^{-1}  \nonumber\\
        & =
        \left(
            \begin{pmatrix}
                \Hh & \mathbf{0} \\
                \B' & \D'
            \end{pmatrix}
            \begin{pmatrix}
                \Hh & \B \\
                \mathbf{0} & \D
            \end{pmatrix} +
            \begin{pNiceArray}{cc|c}
                n_C & 0 & \onevec'\Xc \\
                0 & n_T & \onevec'\Xt \\
                \hline
                \Xc'\onevec & \Xt'\onevec & \Xc'\Xc + \Xt'\Xt
            \end{pNiceArray}
        \right)^{-1}  \nonumber\\
        & =
        \begin{pmatrix}
            \Hh^2 & \Hh\B \\
            \B'\Hh & \B'\B + \D'\D
        \end{pmatrix} +
        \begin{pNiceArray}{cc|c}
            n_C & 0 & n_C \cdot \xbc \\
            0 & n_T & n_T \cdot \xbt \\
            \hline
            n_C \cdot \xbc' & n_T \cdot \xbt' & \X'\X
        \end{pNiceArray} \\
        & =
        \begin{pNiceArray}{cc|c}
            h_1^2 + n_C & 0 & h_1 \bvec_1 + n_C \cdot \xbc \\
            0 & h_2^2 + n_T & h_2 \bvec_2 + n_T \cdot \xbt \\
            \hline
            (h_1 \bvec_1 + n_C \cdot \xbc)' & (h_2 \bvec_2 + n_T \cdot \xbt)' & \B'\B + \D'\D + \X'\X
        \end{pNiceArray}^{-1}. \nonumber
    \end{align}
    Note that $\xbc=\frac{\onevec'\Xc}{n_C}$ and $\xbt=\frac{\onevec'\Xt}{n_T}$.

    Define
    \begin{equation} 
        \Ll =
        \begin{pmatrix}
            h_1^2 + n_C & 0 \\
            0 & h_2^2 + n_T
        \end{pmatrix}^{-1}
        = \begin{pmatrix}
            \frac{1}{s_C} & 0 \\
            0 & \frac{1}{s_T}
        \end{pmatrix}.
    \end{equation}

    The $2 \times 2$ upper left block of $\V_1 = (\Q'\Q + \Z'\Z)^{-1}$, denoted $\R$, is given by
    \begin{align} 
        \R  \nonumber\\
        = & \Big[
            \begin{pmatrix}
                h_1^2 + n_C & 0 \\
                0 & h_2^2 + n_T
            \end{pmatrix}
            - \begin{pmatrix}
                h_1 \bvec_1 + n_C \cdot \xbc \\
                h_2 \bvec_2 + n_T \cdot \xbt
            \end{pmatrix}
            \bigl(
                \B'\B + \X'\X + \D'\D
            \bigr)^{-1}
            \begin{pmatrix}
                h_1 \bvec_1 + n_C \cdot \xbc \\
                h_2 \bvec_2 + n_T \cdot \xbt
            \end{pmatrix}'
        \Big]^{-1}  \nonumber\\
        = & \Big[
            \Ll^{-1}
            - \begin{pmatrix}
                s_C \cdot \gbc \\
                s_T \cdot \gbt
            \end{pmatrix}
            \bigl(
                \B'\B + \X'\X + \D'\D
            \bigr)^{-1}
            \begin{pmatrix}
                s_C \cdot \gbc \\
                s_T \cdot \gbt
            \end{pmatrix}'
        \Big]^{-1}.
        \label{R_before_inversion}
    \end{align}

    By the Woodbury matrix identity \cite[p.~108]{bernstein_matrix_2009},
    \begin{align} 
        \R \nonumber\\
        = & \Ll  - \Ll
        \begin{pmatrix}
            s_C \cdot \gbc \\
            s_T \cdot \gbt
        \end{pmatrix}
        \Big[
            - (\B'\B + \X'\X + \D'\D)
            + \begin{pmatrix}
                s_C \cdot \gbc \\
                s_T \cdot \gbt
            \end{pmatrix}'
            \Ll
            \begin{pmatrix}
                s_C \cdot \gbc \\
                s_T \cdot \gbt
            \end{pmatrix}
        \Big]^{-1}
        \begin{pmatrix}
            s_C \cdot \gbc \\
            s_T \cdot \gbt
        \end{pmatrix}'
        \Ll \nonumber\\
        = & \Ll
        + \begin{pmatrix}
            \gbc \\
            \gbt
        \end{pmatrix}
        \Big[
            (\B'\B + \X'\X + \D'\D)
            - \begin{pmatrix}
                \sqrt{s_C} \gbc \\
                \sqrt{s_T} \gbt
            \end{pmatrix}'
            \begin{pmatrix}
                \sqrt{s_C} \gbc \\
                \sqrt{s_T} \gbt
            \end{pmatrix}
        \Big]^{-1}
        \begin{pmatrix}
            \gbc \\
            \gbt
        \end{pmatrix}' \\
        = & \Ll
        + \begin{pmatrix}
            \gbc \\
            \gbt
        \end{pmatrix}
        \varphi^{-1}
        \begin{pmatrix}
            \gbc \\
            \gbt
        \end{pmatrix}' \nonumber
    \end{align}
    where $\varphi = \Big[
        (\B'\B + \X'\X + \D'\D)
        - \begin{pmatrix}
            \sqrt{s_C} \gbc \\
            \sqrt{s_T} \gbt
        \end{pmatrix}'
        \begin{pmatrix}
            \sqrt{s_C} \gbc \\
            \sqrt{s_T} \gbt
        \end{pmatrix}
    \Big]$.

    Applying the contrast matrix $\A$, we arrive at
    \begin{align} 
        \A \V_1 \A'
        & =
        \begin{pmatrix}
            -1 & 1
        \end{pmatrix}
        \R
        \begin{pmatrix}
            -1 \\
            1
        \end{pmatrix} \nonumber\\
        & = \begin{pmatrix}
            -1 & 1
        \end{pmatrix}
        \left(
            \Ll
            + \begin{pmatrix}
                \gbc \\
                \gbt
            \end{pmatrix}
            \varphi^{-1}
            \begin{pmatrix}
                \gbc \\
                \gbt
            \end{pmatrix}'
        \right)
        \begin{pmatrix}
            -1 \\
            1
        \end{pmatrix} \label{variance_dependent_on_phi}\\
        & = \frac{1}{s_C} + \frac{1}{s_T}
        + (\gbt - \gbc)
        \varphi^{-1}
        (\gbt - \gbc)'. \nonumber
    \end{align}

    Using the following identities:
    \begin{align} 
        \gb - \gbc = \frac{(s_C - s) \gbc + s_T \gbt}{s} = \frac{s_T}{s} (\gbt - \gbc) \\
        \gb - \gbt = \frac{s_C \gbc + (s_T - s) \gbt}{s} = \frac{s_C}{s} (\gbc - \gbt),
    \end{align}
    we obtain
    \begin{align}
        \begin{pmatrix}
            \sqrt{s_C} \gbc \\
            \sqrt{s_T} \gbt
        \end{pmatrix}'
        \begin{pmatrix}
            \sqrt{s_C} \gbc \\
            \sqrt{s_T} \gbt
        \end{pmatrix}
        & = s_C \gbc'\gbc + s_T \gbt'\gbt \nonumber\\
        & = s \gb'\gb - s \gb'\gb + s_C \gbc'\gbc + s_T \gbt'\gbt \nonumber\\
        & = s \gb'\gb - s \gb'\left( \frac{s_C}{s} \gbc + \frac{s_T}{s} \gbt \right) + s_C \gbc'\gbc + s_T \gbt'\gbt \nonumber\\
        & = s \gb'\gb - s_C \gb'\gbc - s_T \gb'\gbt + s_C \gbc'\gbc + s_T \gbt'\gbt \\
        & = s \gb'\gb - s_C (\gb - \gbc)'\gbc - s_T (\gb - \gbt)'\gbt \nonumber\\
        & = s \gb'\gb - s_C \frac{s_T}{s} (\gbt - \gbc)'\gbc - s_T \frac{s_C}{s} (\gbc - \gbt)'\gbt \nonumber\\
        & = s \gb'\gb + \frac{s_C s_T}{s} (\gbt - \gbc)'(\gbt - \gbc) \nonumber
    \end{align}
    and thus
    \begin{align} 
        \varphi = \B'\B + \X'\X - s \gb'\gb - \frac{s_C s_T}{s} (\gbt - \gbc)'(\gbt - \gbc) + \D'\D.
    \end{align}

    Let $\uu = (\gbt - \gbc)'$ be the difference-in-means $p \times 1$ column vector and
    \begin{align} 
        \K = - \frac{s}{s_C s_T}( \B'\B + \X'\X - s \gb'\gb + \D'\D).
    \end{align}
    Then
    \begin{align}
        \varphi^{-1} = -\frac{s}{s_C s_T}(K + \uu\uu')^{-1} = -\frac{s}{s_C s_T}\left(\K^{-1} - \frac{\K^{-1}\uu\uu'\K^{-1}}{1+\uu'\K^{-1}\uu}\right)
    \end{align}
    by the Sherman-Morrison formula \cite[p.~141]{bernstein_matrix_2009}.

    Per (\ref{variance_dependent_on_phi}),
    \begin{align} 
        \A \V_1 \A'
        & = \frac{1}{s_C} + \frac{1}{s_T}
        + (\gbt - \gbc)
        \varphi^{-1}
        (\gbt - \gbc)' \nonumber\\
        & = \frac{s_T}{s_C s_T} + \frac{s_C}{s_C s_T}
        - \frac{s}{s_C s_T}
        \uu'
        \left(
            \K^{-1}
            - \frac{\K^{-1}\uu\uu'\K^{-1}}{1+\uu'\K^{-1}\uu}
        \right)
        \uu \nonumber\\
        & = \frac{s}{s_C s_T}
        - \frac{s}{s_C s_T}
        \left(
            \uu'\K^{-1}\uu
            - \frac{\uu'\K^{-1}\uu\uu'\K^{-1}\uu}{1+\uu'\K^{-1}\uu}
        \right) \nonumber\\
        & = \frac{s}{s_C s_T}
        \left(
            \frac{1+\uu'\K^{-1}\uu}{1+\uu'\K^{-1}\uu}
            - \frac{\uu'\K^{-1}\uu + (\uu'\K^{-1}\uu)^2}{1+\uu'\K^{-1}\uu}
            + \frac{(\uu'\K^{-1}\uu)^2}{1+\uu'\K^{-1}\uu}
        \right) \\
        & = \frac{s}{s_C s_T}
        \left(
            \frac{1}{1+\uu'\K^{-1}\uu}
        \right) \nonumber\\
        & = \frac{s}{s_C s_T}
        \left(
            \frac{1}{1 + (\gbt - \gbc)(-\frac{s_C s_T}{s} (\B'\B + \X'\X - s \gb'\gb + \D'\D)^{-1})(\gbt - \gbc)'}
        \right) \nonumber\\
        & =
        \frac{(\frac{s}{s_C s_T})^2}{\frac{s}{s_C s_T} - (\gbt - \gbc)(\B'\B + \X'\X - s \gb'\gb + \D'\D)^{-1}(\gbt - \gbc)'}. \nonumber
    \end{align}

    We finish the proof by remembering from Lemma~\ref{lem:bayesian_risk} that
    \begin{align}
        r_{\X}(\w) & = \A \V_1 \A' E[\sigma^2] \nonumber\\
        & = \frac{(\frac{s}{s_C s_T})^2}{\frac{s}{s_C s_T} - (\gbt - \gbc)(\B'\B + \X'\X - s \gb'\gb + \D'\D)^{-1}(\gbt - \gbc)'} E[\sigma^2].
    \end{align}
\end{proof}

\begin{proof}[Proof of Corollary~\ref{col:risk_corollary_with_prior}]
    Note that $\gb$ is the mean vector of $\G$, and remember that $\B = \begin{pmatrix}
        \bvec_1 \\
        \bvec_2
    \end{pmatrix}$, so
    \begin{align}
        \B'\B
        = \bvec_1'\bvec_1 + \bvec_2'\bvec_2
        = \sum_{i=1}^{h_1^2} (\frac{\bvec_1}{h_1})'(\frac{\bvec_1}{h_1}) + \sum_{i=1}^{h_2^2} (\frac{\bvec_2}{h_2})'(\frac{\bvec_2}{h_2})
        = \begin{pmatrix}
            \bvec_1/h_1 \\
            \vdots \\
            \bvec_1/h_1 \\
            \bvec_2/h_2 \\
            \vdots \\
            \bvec_2/h_2
        \end{pmatrix}'
        \begin{pmatrix}
            \bvec_1/h_1 \\
            \vdots \\
            \bvec_1/h_1 \\
            \bvec_2/h_2 \\
            \vdots \\
            \bvec_2/h_2
        \end{pmatrix}.
    \end{align}
    It follows that $\X'\X + \B'\B = \G'\G$ and $\B'\B + \X'\X - s \gb'\gb = \Sca(\G)$.
\end{proof}

\begin{proof}[Proof of Corollary~\ref{col:risk_flat_prior}]
    Under the flat prior limit $\V_0^{-1} \rightarrow \mathbf{0}_{(p+2) \times (p+2)}$, equivalent to $f(\zeta \mid \sigma^2) \; \alpha \; 1$, $\V_1 = (\Z'\Z)^{-1}$. We apply Theorem~\ref{thm:risk_theorem_with_prior} noting $s_C=n_C$, $s_T=n_T$, and $\B'\B=\D'\D=\mathbf{0}_{(p+2) \times (p+2)}$.
\end{proof}

\begin{proof}[Proof of Theorem~\ref{thm:sequential_risk}]
    The risk of decision $\w_2$ given new covariate matrix $\U$ and the results of previous observations $\X, \y_1, \w_1$ is given by
    \begin{align} 
        & r_{\U, \X, \y_1, \w_1}(\w_2) \nonumber\\
        = & E[Var(\tau \mid \U, \y_2, \w_2, \X, \y_1, \w_1, \sigma^2) \mid \U, \X, \y_1, \w_1] \nonumber\\
        = & \int_{\sigma^2} \int_{\y_2} Var(\tau \mid \U, \y_2, \w_2, \X, \y_1, \w_1, \sigma^2) f(\y_2 \mid \U, \w_2, \X, \y_1, \w_1, \sigma^2) d\y_2 \nonumber\\
        & \qquad \qquad \cdot f(\sigma^2 \mid \U, \w_2, \X, \y_1, \w_1) d\sigma^2 \nonumber\\
        = & \int_{\zeta} \int_{\sigma^2} \int_{\y_2} Var(\tau \mid \U, \y_2, \w_2, \X, \y_1, \w_1, \sigma^2) f(\y_2 \mid \U, \w_2, \X, \y_1, \w_1, \zeta, \sigma^2) d\y_2 \nonumber\\
        & \qquad \qquad \cdot f(\zeta, \sigma^2 \mid \U, \w_2, \X, \y_1, \w_1) d\sigma^2 d\zeta \nonumber\\
        = & \int_{\zeta} \int_{\sigma^2} \int_{\y_2} \A Var(\zeta \mid \U, \y_2, \w_2, \X, \y_1, \w_1, \sigma^2) \A' \; f(\y_2 \mid \U, \w_2, \zeta, \sigma^2) d\y_2 \nonumber\\
        & \qquad \qquad \cdot f(\zeta, \sigma^2 \mid \X, \y_1, \w_1) d\sigma^2 d\zeta \nonumber\\
        = & \int_{\zeta} \int_{\sigma^2} \int_{\y_2} \A \V_2 \A' \sigma^2 \; f(\y_2 \mid \U, \w_2, \zeta, \sigma^2) d\y_2 \; f(\zeta, \sigma^2 \mid \X, \y_1, \w_1) d\sigma^2 d\zeta \label{risk_sequential_1}\\
        = & \A \V_2 \A' \int_{\zeta} \int_{\sigma^2} \int_{\y_2} f(\y_2 \mid \U, \w_2, \zeta, \sigma^2) d\y_2 \; \sigma^2 \; f(\zeta, \sigma^2 \mid \X, \y_1, \w_1) d\sigma^2 d\zeta \nonumber\\
        = & \A \V_2 \A' E[\sigma^2 \mid \X, \y_1, \w_1] \nonumber
    \end{align}
    where $f(\y_2 \mid \U, \w_2, \X, \y_1, \w_1, \zeta, \sigma^2) = f(\y_2 \mid \U, \w_2, \zeta, \sigma^2)$ by conditional independence of the $y_i$ given $\mathbf{x}_i, \w_i, \zeta, \sigma^2$ and $f(\zeta, \sigma^2 \mid \U, \w_2, \X, \y_1, \w_1) = f(\zeta, \sigma^2 \mid \X, \y_1, \w_1)$ by
    \begin{align}
        & \; f(\zeta, \sigma^2 \mid \U, \w_2, \X, \y_1, \w_1) \nonumber\\
        = & \int_{\psi} f(\zeta, \sigma^2, \psi \mid \U, \w_2, \X, \y_1, \w_1) d\psi \nonumber\\
        \alpha & \int_{\psi} f(\U, \w_2, \X, \y_1, \w_1 \mid \zeta, \sigma^2, \psi) f(\zeta, \sigma^2, \psi) d\psi \nonumber\\
        = & \int_{\psi} f(\X, \y_1, \w_1 \mid \zeta, \sigma^2, \psi) f(\U \mid \zeta, \sigma^2, \psi) f(\zeta, \sigma^2) f(\psi) d\psi \\
        = & \int_{\psi} f(\y_1 \mid \X, \w_1, \zeta, \sigma^2) f(\X \mid \psi) f(\U \mid \psi) f(\zeta, \sigma^2) f(\psi) d\psi \nonumber\\
        = & \; f(\y_1 \mid \X, \w_1, \zeta, \sigma^2) f(\zeta, \sigma^2) \int_{\psi} f(\X \mid \psi) f(\U \mid \psi) f(\psi) d\psi \nonumber\\
        \alpha & \; f(\zeta, \sigma^2 \mid \X, \y_1, \w_1). \nonumber
    \end{align}

    $\w_1$ and $\w_2$ are taken to be deterministic functions of $\X$ and $(\X, \U)$, respectively. A consequence of (\ref{risk_sequential_1}) is that $r_{\U, \X, \y_1, \w_1}(\w_2) = \A \V_2 \A' E[\sigma^2 \mid \X, \y_1, \w_1]$ depends on $\w_2$ only through $\A \V_2 \A'$.

    Under the prior distribution for the first allocation, $\V_1 = (\Q'\Q + \Z'\Z)^{-1}$, where $\Z = \begin{pNiceArray}{cc|c}
        \onevec & \zerovec & \Xc \\
        \zerovec & \onevec & \Xt
    \end{pNiceArray}$. Denote $\J = \begin{pNiceArray}{cc|c}
        \onevec & \zerovec & \Uc \\
        \zerovec & \onevec & \Ut
    \end{pNiceArray}$. All vectors $\zerovec$ and $\onevec$ have the conforming size. Then
    \begin{align} 
        \V_2
        = & (\V_1^{-1} + \J'\J)^{-1} = (\Q'\Q + \Z'\Z + \J'\J)^{-1} \nonumber\\
        = & \left(
            \begin{pmatrix}
                \Hh^2 & \Hh\B \\
                \B'\Hh & \B'\B + \D'\D
            \end{pmatrix} \right. \nonumber\\
            & \left. \;\;\; +
            \begin{pNiceArray}{cc|c}
                n_C & 0 & \onevec'\Xc \\
                0 & n_T & \onevec'\Xt \\
                \hline
                \Xc'\onevec & \Xt'\onevec & \Xc'\Xc + \Xt'\Xt
            \end{pNiceArray} +
            \begin{pNiceArray}{cc|c}
                m_C & 0 & \onevec'\Uc \\
                0 & m_T & \onevec'\Ut \\
                \hline
                \Uc'\onevec & \Ut'\onevec & \Uc'\Uc + \Ut'\Ut
            \end{pNiceArray}
        \right)^{-1} \nonumber\\
        = &
        \left(
            \begin{pNiceArray}{cc|c}
                h_1^2 & 0 & h_1 \bvec_1 \\
                0 & h_2^2 & h_2 \bvec_2 \\
                \hline
                h_1 \bvec_1' & h_2 \bvec_2' & \B'\B + \D'\D
            \end{pNiceArray} \right. \nonumber\\
            & \left. \;\;\; +
            \begin{pNiceArray}{cc|c}
                n_C + m_C & 0 & \onevec'\Xc + \onevec'\Uc \\
                0 & n_T + m_T & \onevec'\Xt + \onevec'\Ut \\
                \hline
                \Xc'\onevec + \Uc'\onevec & \Xt'\onevec + \Ut'\onevec & \Xc'\Xc + \Xt'\Xt + \Uc'\Uc + \Ut'\Ut
            \end{pNiceArray}
        \right)^{-1} \nonumber\\
        = &
        \left(
            \begin{pNiceArray}{cc|c}
                h_1^2 & 0 & h_1 \bvec_1 \\
                0 & h_2^2 & h_2 \bvec_2 \\
                \hline
                h_1 \bvec_1' & h_2 \bvec_2' & \B'\B + \D'\D
            \end{pNiceArray} \right. \nonumber\\
            & \left. \;\;\; +
            \begin{pNiceArray}{cc|c}
                l_C & 0 & \onevec'\Xc + \onevec'\Uc \\
                0 & l_T & \onevec'\Xt + \onevec'\Ut \\
                \hline
                \Xc'\onevec + \Uc'\onevec & \Xt'\onevec + \Ut'\onevec & \Xc'\Xc + \Uc'\Uc + \Xt'\Xt + \Ut'\Ut
            \end{pNiceArray}
        \right)^{-1} \nonumber\\
        = &
        \left(
            \begin{pNiceArray}{cc|c}
                h_1^2 & 0 & h_1 \bvec_1 \\
                0 & h_2^2 & h_2 \bvec_2 \\
                \hline
                h_1 \bvec_1' & h_2 \bvec_2' & \B'\B + \D'\D
            \end{pNiceArray} +
            \begin{pNiceArray}{cc|c}
                l_C & 0 & l_C \obc \\
                0 & l_T & l_T \obt \\
                \hline
                l_C \obc' & l_T \obt' & \Oo'\Oo
            \end{pNiceArray}
        \right)^{-1} \\
        = &
        \begin{pNiceArray}{cc|c}
            h_1^2 + l_C & 0 & h_1 \bvec_1 + l_C \obc \\
            0 & h_2^2 + l_T & h_2 \bvec_2 + l_T \obt \\
            \hline
            (h_1 \bvec_1 + l_C \obc)' & (h_2 \bvec_2 + l_T \obt)' & \B'\B + \D'\D + \Oo'\Oo
        \end{pNiceArray}^{-1} \nonumber\\
        = &
        \begin{pNiceArray}{cc|c}
            s_C & 0 & s_C \gbc \\
            0 & s_T & s_T \gbt \\
            \hline
            s_C \gbc' & s_T \gbt' & \B'\B + \D'\D + \Oo'\Oo
        \end{pNiceArray}^{-1}. \nonumber
    \end{align}

    Applying the contrast matrix $\A$ to $\V_2$, all terms except for the upper left $2 \times 2$ submatrix are zeroed out, resulting in the following by block matrix inversion:
    \begin{align}
        & \A \V_2 \A' \nonumber\\
        = & \begin{pmatrix}
            -1 & 1
        \end{pmatrix}
        \Big[
            \begin{pmatrix}
                s_C & 0 \\
                0 & s_T
            \end{pmatrix}
            - \begin{pmatrix}
                s_C \gbc \\
                s_T \gbt
            \end{pmatrix}
            \bigl(
                \B'\B + \D'\D + \Oo'\Oo
            \bigr)^{-1}
            \begin{pmatrix}
                s_C \gbc \\
                s_T \gbt
            \end{pmatrix}'
        \Big]^{-1}
        \begin{pmatrix}
            -1 \\
            1
        \end{pmatrix}.
    \end{align}

    Note the inner matrix is the same as in (\ref{R_before_inversion}) with $\Oo$ instead of $\X$. Consequently, all further results in the earlier proof follow.
\end{proof}

\begin{proof}[Proof of Corollary~\ref{col:risk_mahalanobis}]
    Note that $\frac{\Sca(\X)}{n-1} = Cov(\X)$. The result follows from inserting $\Sca(\X) = (n-1) Cov(\X)$ into $r_{\X}(\w) = \frac{(\frac{n}{n_C n_T})^2}{\frac{n}{n_C n_T} - (\xbt - \xbc)\Sca(\X)^{-1}(\xbt - \xbc)'} E[\sigma^2]$ obtained in Corollary~\ref{col:risk_flat_prior} and simplifying the fraction.
\end{proof}

\begin{proof}[Proof of Theorem~\ref{thm:equal_split_condition}]
    Let $\w_1$ and $\w_2$ be two allocations and define $p_1=\frac{n_{T1}}{n}=\frac{\sum \w_1}{n}$ and $p_2=\frac{n_{T2}}{n}=\frac{\sum \w_2}{n}$. Consequently, $\frac{n}{n_{C1}n_{T1}}=\frac{1}{np_1(1-p_1)}$ and $\frac{n}{n_{C2}n_{T2}}=\frac{1}{np_2(1-p_2)}$. Note that for an allocation $\w$ the difference in covariate means vector is given by
    \begin{align}
    \xbt-\xbc
    =& \frac{\onevec'\Xt}{n_T} - \frac{\onevec'\Xc}{n_C} \nonumber\\
    =& \frac{\w'\X}{n_T} - \frac{(\onevec-\w)'\X}{n_C} \nonumber\\
    =& \frac{1}{n_Cn_T}n_C\w'\X -\frac{1}{n_Cn_T}n_T(\onevec-\w)'\X \\
    =& \frac{1}{n_Cn_T}\left(n_C\w + n_T\w - n_T\onevec \right)'\X \nonumber\\
    =& \frac{n}{n_Cn_T}\left(\w - \frac{n_T}{n}\onevec \right)'\X. \nonumber
    \end{align}

    We will assume, without loss of generality, $\X$ has mean 0 in all columns. This can be achieved by centering $\X$, as both the difference in means vector and $\Sca(\X)$ are invariant to translations.
    
    Then:
    \begin{align}
        & \;\;\; r_{\X}(\w_1) \leq r_{\X}(\w_2) \nonumber\\
        \iff & \frac{(\frac{n}{n_{C1} n_{T1}})^2}{\frac{n}{n_{C1} n_{T1}} - (\xbt_1 - \xbc_1)\Sca(\X)^{-1}(\xbt_1 - \xbc_1)'} E[\sigma^2] \nonumber\\
        & \qquad \qquad \leq \qquad \qquad \frac{(\frac{n}{n_{C2} n_{T2}})^2}{\frac{n}{n_{C2} n_{T2}} - (\xbt_2 - \xbc_2)\Sca(\X)^{-1}(\xbt_2 - \xbc_2)'} E[\sigma^2] \nonumber\\
        \iff & \frac{(\frac{n}{n_{C1} n_{T1}})^2}{\frac{n}{n_{C1} n_{T1}} - (\frac{n}{n_{C1} n_{T1}})^2\left(\w_1 - p_1\onevec \right)'\X\Sca(\X)^{-1}\X'\left(\w_1 - p_1\onevec \right)} \nonumber\\
        & \qquad \qquad \leq \qquad \qquad \frac{(\frac{n}{n_{C2} n_{T2}})^2}{\frac{n}{n_{C2} n_{T2}} - (\frac{n}{n_{C2} n_{T2}})^2\left(\w_2 - p_2\onevec \right)'\X\Sca(\X)^{-1}\X'\left(\w_2 - p_2\onevec \right)}  \nonumber\\
        \iff & \frac{1}{(\frac{n}{n_{C1} n_{T1}})^{-1} - \left(\w_1 - p_1\onevec \right)'\X\Sca(\X)^{-1}\X'\left(\w_1 - p_1\onevec \right)} \\
        & \qquad \qquad \leq \qquad \qquad \frac{1}{(\frac{n}{n_{C2} n_{T2}})^{-1} - \left(\w_2 - p_2\onevec \right)'\X\Sca(\X)^{-1}\X'\left(\w_2 - p_2\onevec \right)}  \nonumber\\
        \iff & np_1(1-p_1) - \left(\w_1 - p_1\onevec \right)'\X\Sca(\X)^{-1}\X'\left(\w_1 - p_1\onevec \right) \nonumber\\
        & \qquad \qquad \geq \qquad \qquad np_2(1-p_2) - \left(\w_2 - p_2\onevec \right)'\X\Sca(\X)^{-1}\X'\left(\w_2 - p_2\onevec \right)  \nonumber\\
        \iff & np_1(1-p_1) - np_2(1-p_2) \nonumber\\
        & \geq (\left(\w_1 - p_1\onevec \right)'\X\Sca(\X)^{-1}\X'\left(\w_1 - p_1\onevec \right)) - \left(\w_2 - p_2\onevec \right)'\X\Sca(\X)^{-1}\X'\left(\w_2 - p_2\onevec \right). \nonumber
    \end{align}

    As established before, we can assume centered $\X$, and so
    \begin{equation}
        \X\Sca(\X)^{-1}\X' = \X(\X'\X - n\zerovec'\zerovec)^{-1}\X'=\X(\X'\X)^{-1}\X'=\Hh(\X).
    \end{equation}

    As $\left(\w_2 - p_2\onevec \right)'\Hh(\X)\left(\w_2 - p_2\onevec \right)$ is a quadratic form of a positive-semidefinite matrix, it is lower bounded by 0. Thus, 
    \begin{equation}
        (\left(\w_1 - p_1\onevec \right)'\Hh(\X)\left(\w_1 - p_1\onevec \right)) \leq np_1(1-p_1) - np_2(1-p_2)
    \end{equation}
    for any $p_2$ implies $r_{\X}(\w_1) \leq r_{\X}(\w_2)$ for any $\w_2$ such that $p_2 \neq p_1$. Note that, due to the positive-semi-definiteness of $\left(\w_1 - p_1\onevec \right)'\Hh(\X)\left(\w_1 - p_1\onevec \right)$, the condition can only be achieved when $p_1$ is closer to $\frac{1}{2}$ than $p_2$ is, as $np(1-p)$ forms a parabola in $p$ which is maximized by $p=\frac{1}{2}$ and the right-hand side must be positive.

    Consider the case of $p_1=\frac{1}{2}$, when $n_{C1}=n_{T1}=\frac{n}{2}$ and $np_1(1-p_1)=\frac{n}{4}$, and, without loss of generality, let $n_{T2}=\frac{n}{2}+k$, $k \geq 1$. Then 
    \begin{align}
        np_2(1-p_2)
        =& n\frac{\frac{n}{2}+k}{n}\frac{\frac{n}{2}-k}{n} \nonumber\\
        =& n(\frac{1}{2}+\frac{k}{n})(\frac{1}{2}-\frac{k}{n}) \\
        =& n(\frac{1}{4}-\frac{k^2}{n^2}) \nonumber\\
        =& \frac{n}{4}-\frac{k^2}{n}. \nonumber
    \end{align}
    The right-hand side of the condition becomes
    \begin{equation}
        np_1(1-p_1) - np_2(1-p_2) = \frac{k^2}{n} \geq \frac{1}{n}
    \end{equation}
    for any $k \geq 1$. It follows that existence of $\w_1$ with $p_1=\frac{1}{2}$ such that
    \begin{equation}
        \left(\w_1 - \frac{1}{2}\onevec \right)'\Hh(\X)\left(\w_1 - \frac{1}{2}\onevec \right) \leq \frac{1}{n}
    \end{equation}
    guarantees $r_{\X}(\w_1) \leq r_{\X}(\w_2)$ for any $\w_2$ such that $p_2 \neq \frac{1}{2}$, and so the optimal allocation $\w$ must have $n_C=n_T=\frac{n}{2}$.
\end{proof}

\end{document}